\documentclass[10pt]{article}
\usepackage{amsmath}
\usepackage{amssymb}
\usepackage{amscd}
\usepackage{latexsym}
\usepackage{theorem}
\usepackage{setspace}
\usepackage{epsfig}
\usepackage{amsfonts}
\usepackage{enumerate}
\usepackage{xypic}
\usepackage{mathrsfs}

\newtheorem{theorem}{Theorem}[section]
\newtheorem{lemma}[theorem]{Lemma}
\newtheorem{proposition}[theorem]{Proposition}
\newtheorem{corollary}[theorem]{Corollary}

{\theorembodyfont{\rmfamily}
\theoremstyle{plain}
\newtheorem{definition}[theorem]{Definition}
\newtheorem{example}[theorem]{Example}

\newtheorem{remark}[theorem]{Remark}
}

\input xy
\xyoption{all}

\renewcommand{\Re}{\operatorname{Re}}

\renewcommand{\dim}{\operatorname{dim}}

\newcommand{\C}{{\mathbb{C}}}
\newcommand{\Z}{{\mathbb{Z}}}
\newcommand{\Q}{{\mathbb{Q}}}

\newcommand{\R}{{\mathbb{R}}}
\newcommand{\Zt}{{\mathbb{Z}_2}}

\newcommand{\g}{\gamma}

\renewcommand{\a}{\alpha}
\renewcommand{\b}{\beta}

\newcommand{\s}{\sigma}
\newcommand{\eps}{\epsilon}

\newcommand{\subs}{\subseteq}

\newcommand{\impl}{\Rightarrow}

\renewcommand{\(}{\left(}
\renewcommand{\)}{\right)}

\renewcommand{\>}{\right>}

\newcommand{\hs}{\hspace{3pt}}

\newcommand{\Sp}{S^+}
\newcommand{\Sm}{S^-}
\newcommand{\A}{\mathcal{A}}

\newcommand{\CA}{\mathcal{C}(\A)}
\newcommand{\MA}{\mathcal{M}(\A)}

\newcommand{\Hip}{H_i^+}
\newcommand{\Him}{H_i^-}

\newcommand{\vgz}{VG(\A;\Z)}

\newcommand{\dA}{d\A}
\newcommand{\Lda}{\Line(\dA)}
\newcommand{\Mda}{\mathcal{M}(\dA)}
\newcommand{\cs}{\C^*}
\newcommand{\Kcs}{\Line(\cs)}
\newcommand{\LA}{\Line(\A)}
\newcommand{\LAS}{\Line(\A_S)}
\newcommand{\PA}{P(\A)}
\newcommand{\IA}{\mathcal{I}_{\A}}
\newcommand{\oi}{\omega_i}
\newcommand{\etai}{\eta_i}
\newcommand{\etaj}{\eta_j}
\newcommand{\etaz}{\eta_0}
\newcommand{\on}{\{1,\ldots,n\}}
\newcommand{\codim}{\operatorname{codim}}

\newcommand{\Line}{Line}
\newcommand{\MAS}{\mathcal{M}(\A_S)}
\newcommand{\Cd}{\C^d}
\newcommand{\D}{\Delta}
\newcommand{\BA}{B(\A)}

\newcommand{\Capr}{\mathcal{C}(\A'')}
\newcommand{\Lca}{\Line\big(\CA\big)}
\newcommand{\kg}{K_{\! G}}
\newcommand{\kog}{K\! O_{\! G}}
\newcommand{\kot}{K\! O_{\Zt}}
\newcommand{\kotx}{\kot\!(X)}

\newcommand{\kotxq}{\kot\!(X)_{\Q}}
\newcommand{\kt}{K_{\Zt}}
\newcommand{\ktx}{\kt\!(X)}
\newcommand{\kx}{K(X)}
\newcommand{\kxq}{K(X)_{\Q}}
\newcommand{\ktxq}{\kt\!(X)_{\Q}}
\newcommand{\half}{\frac{1}{2}}
\newcommand{\fo}{f\!o}

\newcommand{\qed}{\hfill \mbox{$\Box$}\medskip\newline}
\newenvironment{proof}{\noindent {\bf Proof:}}{\qed \par}
\newenvironment{proofoverq}{\noindent {\bf Proof of Proposition \ref{overq}:}}{\qed \par}

\begin{document}
\begin{spacing}{1.2}

\noindent
{\LARGE \bf Hyperplane arrangements and K-theory}
\bigskip\\
{\bf Nicholas Proudfoot}\footnote{Partially supported
by the Clay Mathematics Institute Liftoff Program} \\
Department of Mathematics, University of California,
Berkeley, CA 94720
\bigskip
{\small
\begin{quote}
\noindent {\em Abstract.} We study the $\Zt$-equivariant K-theory
of $\MA$, where $\MA$ is the complement of the complexification of a
real hyperplane arrangement, and $\Zt$ acts on $\MA$ by complex
conjugation.  We compute the rational equivariant K- and KO-rings
of $\MA$, and we give two different combinatorial descriptions
of a subring $\LA$ of the integral equivariant KO-ring, where
$\LA$ is defined to be the subring generated by equivariant line
bundles.
\end{quote}
}

\begin{section}{Introduction}\label{intro}
Let $\A$ be an arrangement of $n$ hyperplanes in $\C^d$, and let
$\MA$ denote the complement of $\A$ in $\C^d$. It is a fundamental
problem in the study of hyperplane arrangements to investigate the
extent to which the topology of $\MA$ is determined by the
combinatorics (more precisely the pointed matroid) of $\A$. Perhaps
the first major theorem in the subject is the celebrated result of
Orlik and Solomon \cite{OS}, in which the cohomology ring of $\MA$
is shown to have a combinatorial presentation in terms of the
pointed matroid.  Our goal is to give a combinatorial description of
the K-theory of $\MA$.

We will work only with
hyperplane arrangements which are defined over the real numbers.
Though restrictive, this hypothesis allows for more subtle constructions
in both combinatorics and topology.
Let $\A = \{H_1,\ldots,H_n\},$ where $H_i$ is the zero set of an affine
linear map $\oi:\R^d\to\R$, and let $H_i^{\pm}=\oi^{-1}(\R^{\pm})$ be
the corresponding {\em open} half-spaces in $\R^d$.
On the combinatorial side, a real hyperplane arrangement determines a
{\em pointed oriented matroid} \cite{BLSWZ}.
The pointed oriented matroid
of $\A$ is characterized by two types of combinatorial data:
\begin{enumerate}
\item which subsets $S\subs\on$ have the property that $\bigcap_{i\in S}H_i$
is nonempty with codimension less than $|S|$, and
\item which pairs of subsets $S^+,S^-\subs\{1,\ldots,n\}$ have the property that
$\bigcap_{i\in S^+}H_i^+\cap\bigcap_{j\in S^-}H_j^- = \varnothing$.
\end{enumerate}
On the topological side, the complement $\MA$ of the complexified arrangement
carries an action of $\Zt$, given by complex conjugation.  This allows us
to consider not only the ordinary algebraic invariants, but their $\Zt$-equivariant
analogues as well.  The equivariant fundamental group and the equivariant
cohomology ring have been studied in \cite{Hu} and \cite{Pr}, respectively.
In \cite{Pr}, we extend a theorem of Salvetti \cite{Sa}
to show that the pointed oriented matroid determines
the equivariant homotopy type of $\MA$, hence the ``extra'' combinatorics
and ``extra'' topology arising from the real structure on $\A$ go hand in hand.

Our main result is to give two combinatorial descriptions of the
ring $\LA$, which we define to be the subring of the degree zero
equivariant KO-ring $\kot\!\big(\MA\big)$ generated by line bundles.
We first present $\LA$ as a quotient of a polynomial ring, in a
manner similar to our presentation of the equivariant cohomology
ring of $\MA$ in \cite{Pr}.  One important difference is that the
equivariant cohomology ring is only well behaved with coefficients
in $\Zt$, whereas $\LA$ is both interesting and computable over the
integers. We then give a second description of $\LA$ as a subring of
the equivariant KO-ring of the fixed point set $\CA$, the complement
of the real arrangement. Since $\CA$ is a finite disjoint union of
contractible spaces, its equivariant KO-ring is simply a direct sum
of equivariant KO-rings of points.

In Section \ref{EK}, we also compute
the more familiar rings $\kot\!\big(\MA\big)$ and $K_{\Zt}\!\big(\MA\big)$
after tensoring with the rational numbers
(see Proposition \ref{overq}, Remark \ref{koa}, and Corollary \ref{cplx}).
A dimension count reveals that these rings are strictly larger than the
tensor product of $\LA$ with $\Q$; in other words, they are {\em not}
entirely generated by line bundles.
We find that these rings may be described purely in terms of the ordinary cohomology
rings of $\MA$ and $\CA$, thus the only truly new invariants
come from working over $\Z$.

\vspace{\baselineskip}
\noindent {\em Acknowledgments.}
I am grateful to Laurent Bartholdi, Allen Knutson, Andr\'e Henriques,
Greg Landweber,
and Alexander Varchenko for helpful discussions.
\end{section}

\begin{section}{Equivariant K-theory}\label{EK}
Let $X$ be a topological space equipped with an action of a group
$G$. The {\em equivariant K-ring} $\kg(X)$ is defined to be the
Grothendieck ring of $G$-equivariant complex vector bundles on $X$.
More precisely, $\kg(X)$ is additively generated by 
$G$-equivariant complex vector bundles over $X$, modulo
the ideal generated by elements of the form $\sum_{i=0}^m (-1)^iE^i$
for every exact sequence $0\to E^1\to E^2\to\ldots\to E^m\to 0.$ The
multiplicative structure is given by the tensor product, and the
trivial line bundle is the multiplicative identity. Then $\kg$ is a
contravariant functor from $G$-spaces to rings, and is constant on
$G$-equivariant homotopy equivalence classes of $G$-spaces. Since
every $G$-space $X$ maps $G$-equivariantly to a point, $\kg(X)$ is
naturally a module over $\kg(pt)$, the representation ring of $G$.
Similarly, we may define the contravariant functor $\kog$ from
$G$-spaces to rings, which takes a space $X$ to its Grothendieck
ring of $G$-equivariant {\em real} vector bundles, which is a module
over the real representation ring $\kog(pt)$.

We will say that a (not necessarily exact) complex of bundles $0\to
E^1\to E^2\to\ldots\to E^m\to 0$ {\em represents} the element
$\sum_{i=0}^m (-1)^i E^i$ in either real or complex K-theory.
Given two complexes, we may tensor them together and then add up the
diagonals to make a third complex, and the class represented by the
tensor product of the two complexes is equal to the product of the
classes represented by each complex. If $E^{\bullet}$ and
$F^{\bullet}$ are complexes of $G$-equivariant vector bundles, then
the locus of points in $X$ over which the tensor product $(E\otimes
F)^{\bullet}$ fails to be exact is contained in the intersection of
the loci over which $E^{\bullet}$ and $F^{\bullet}$ individually
fail to be exact. In particular, given any two K-theory classes
which may be represented by complexes that fail to be exact on
disjoint sets, their product is equal to zero.  This will be our
principal means of identifying relations in $\kog(X)$ (see Example
\ref{cstar} and Theorem \ref{pl}).

In this paper we will be concerned only with the case $G=\Zt$.  Furthermore, we will
restrict our attention in later sections
to a subring $\Line(X)\subs\kotx$, which we define
to be the subring additively generated by line bundles.
Though not part of a generalized cohomology theory, $\Line$ is a contravariant functor from $\Zt$-spaces to rings,
and $\Line(X)$ is always a module over $\Line(pt)$.
The ring $\Line(pt)$ is additively generated by the unit element $1$, and the element
$N\in\Line(pt)$ representing the unique nontrivial one-dimensional
representation of $\Zt$, subject to the relation $N^2=1$.
We will write $x=1-N$, so that we have
$$\Line(pt) = \Z[x]/x(2-x).$$
For an arbitrary $\Zt$-space $X$, we will abuse notation by writing
$x\in\Line(X)$ to denote the image of $x\in\Line(pt)$.

Real $\Zt$-equivariant line bundles on $X$ are classified by the
equivariant cohomology group $H^1_{\Zt}\!(X;\Zt)$, with the
isomorphism given by the first equivariant Stiefel-Whitney class.
(The completely analogous statement for complex line bundles is
proven in \cite[C.6.3]{GGK} as well as \cite[A.1]{HL}.)
Hence $\Line(X)$ is isomorphic to a
quotient of the group ring $\Z\big[H^1_{\Zt}\!(X;\Zt)\big]$ by
relations that arise when two different sums of equivariant line
bundles are isomorphic. The obvious advantage of working with this
subring is that the group $H^1_{\Zt}\!(X;\Zt)$ is often computable.

Despite the relative intractability of computing the more familiar rings $\kt(X)$
and $\kot(X)$, it is not so hard to compute their rationalizations $\ktxq:=\ktx\otimes\Q$
and $\kotxq:=\kotx\otimes\Q$,
especially in the case where $X$ is the complement of a hyperplane arrangement.
We include this computation here, though it will not be relevant
to the rest of the paper.

Let $\s:X\to X$ be the involution given by the $\Zt$-action. By
definition, a $\Zt$-equivariant vector bundle on $X$ is an ordinary
bundle $E$ along with a choice of isomorphism $E\cong \s^*E$, hence
the image of the forgetful map $\fo:\ktxq\to\kxq$ is equal to the
invariant ring $\kxq^{\s^*}$. Let $x=1-N_{\C}\in\ktxq$, where
$N_{\C}$ is the complexification of the real line bundle $N$ defined
above.

\begin{lemma}\label{forget}
The kernel of the forgetful map $\fo:\ktxq\to\kxq$ is generated by $x$.
\end{lemma}

\begin{proof}
The element $x=1-N$ is clearly contained in the kernel of the forgetful map.
To prove the other containment, we observe the fact that
an equivariant bundle on a free $G$-space carries
the same data as an ordinary bundle on the quotient, hence
the nonequivariant K-ring $\kx$ may be identified with the equivariant
ring $\kt\!(X\times\Zt)$, where $\Zt$ acts diagonally (and therefore freely)
on $X\times\Zt$.  In this picture, the forgetful map gets identified
with the pullback along the projection $\pi:X\times\Zt\to X$.

Consider the pushforward $\pi_*:\kt\!(X\times\Zt)_{\Q}\to\ktxq$,
taking a bundle $E$ on $X\times\Zt$ to $E|_{X\times\{1\}}\oplus E|_{X\times\{-1\}}$.
It is easy to check that the equivariant structure on $E$ defines a natural
equivariant structure on $\pi_*(E)$, and that this pushforward satisfies the
projection formula $\pi_*(\b \cdot\pi^*\a) = \pi_*(\b)\cdot\a$.
Suppose that $\a\in\ker\pi^*$.  Then
$$(2-x)\cdot\a = \pi_*(1)\cdot\a = \pi_*\(\pi^*(\a)\) = 0,$$ hence
$$\a = \half(2-x+x)\a = \half x\a$$ is a multiple of $x$.
\end{proof}

\vspace{-\baselineskip}
\begin{remark}
To prove Lemma \ref{forget} we did not really have to work over the rationals,
we only had to invert $2$.  The analogous statement over the integers is false.
\end{remark}

\begin{proposition}\label{overq}
There is a ring isomorphism
$\ktxq\,\,\cong\,\, H^{2*}\!(X;\Q)^{\s^*}\oplus\,\,\kt(X^{\s})_{\Q}
\big/\!\<x-2\>$.
\end{proposition}

\begin{remark}\label{koa}
Suppose that $X = \MA$ is the complement of the complexification
of a real hyperplane arrangement $\A$, and let $\s$ be the involution given by complex conjugation.
Then the cohomology ring of $X$ is isomorphic to the Orlik-Solomon
algebra of $\A$ \cite{OS}.  The involution $\s^*$ acts by negation
on the generators of the Orlik-Solomon algebra, and therefore the
invariant ring $H^{2*}\!(X;\Q)^{\s^*}$ is simply the even degree
part of the Orlik-Solomon algebra.  Let $\CA = \MA^{\s}$ be the
complement of the real arrangement. This space is a disjoint union
of contractible pieces, hence $\kt(X^{\s})_{\Q}/\!\<x-2\>$ is isomorphic
a product of copies of $\Q = \kt(pt)_{\Q}/\!\<x-2\>$ for each component.
This ring is also known as the Varchenko-Gelfand ring $VG(\A;\Q)$ of
locally constant $\Q$-valued functions on $\CA$.
\end{remark}

\begin{proofoverq}
We begin by considering the map
$$\ktxq\to\ktxq\big/\!\<x\>\oplus\ktxq\big/\!\<x-2\>$$ given by the two projections.
This map is surjective because the generator of the kernel of the first
projection maps to a unit in the second factor.  It is also injective,
because any element of the kernel is annihilated both by $2-x$ and by $x$,
and therefore also by $2$.  By Lemma \ref{forget},
$\ktxq/\!\<x\>$ is isomorphic to $\kxq^{\s^*}$,
and $H^{2*}\!(X;\Q)$ is isomorphic to $\kxq$
via the Chern character.
Hence the first factor of $\ktxq$ is isomorphic to $H^{2*}\!(X;\Q)^{\s^*}$.
The fact that $\ktxq/\!\<x-2\>\cong\kt(X^{\s})/\!\<x-2\>$ is a consequence
of the localization theorem \cite[3.4.1]{AS}.
\end{proofoverq}

\vspace{-\baselineskip}
\begin{corollary}\label{cplx}
There is a ring isomorphism
$\kotxq\,\,\cong\,\, H^{4*}\!(X;\Q)^{\s^*}\oplus\,\,\kot(X^{\s})\big/\!\<x-2\>$.
\end{corollary}

\begin{proof}
The complexification map from $\kotxq$ to $\ktxq$ is injective, and
its image may be identified with the fixed point set under the
involution taking a complex vector bundle to its conjugate \cite[p.
74]{Bo}.  (This involution is not to be confused with the involution
$\s^*$.)  On $H^{2k}(X)$, this involution translates into
multiplication by $(-1)^k$, hence the invariant ring is
$H^{4*}\!(X;\Q)^{\s^*}$.
\end{proof}
\end{section}

\begin{section}{The quotient description of \boldmath$\LA$}\label{la}
Let $\omega_1,\ldots,\omega_n$ be a collection of affine linear
functionals on $\R^d$, and let $\A=\{H_1,\ldots,H_n\}$ be the
associated cooriented hyperplane arrangement.  By the word
cooriented, we mean that we have not only a collection of
hyperplanes, but also a collection of positive open half spaces
$\Hip = \oi^{-1}(\R^+)$, along with their negative counterparts
$\Him=\oi^{-1}(\R^-)$.  Let $\CA = \R^d\smallsetminus\cup_{i=1}^n
H_i$ be the complement of $\A$ in $\R^d$, and let $\MA =
\C^d\smallsetminus\cup_{i=1}^n H_i^{\C}$ be its complexification.
Then $\MA$ carries an action of $\Zt$ given by complex conjugation,
with fixed point set $\CA$. For each $i$, the complexification of
$\oi$ restricts to a map $\MA\to\cs$.  We will abuse notation by
calling this map $\oi$ as well.

The purpose of this section is to give a combinatorial presentation
of the ring $\LA := \Line\big(\MA\big)$.
We begin with the most basic example, where $\A$ consists of
a single point on a line, and therefore $\MA = \cs$.
This example will be fundamental to understanding the general case,
as all line bundles on a general $\MA$ will be constructed
as tensor products of pullbacks of line bundles on $\cs$ along
the maps $\oi:\MA\to\cs$.

\begin{example}\label{cstar}
Let $\A$ consist of one point in $\R$, so that $\MA\cong\cs$. Let
$N$ be the topologically trivial real line bundle on $\cs$ with the
nontrivial $\Zt$-action at every fixed point (the pullback of the
nontrivial $\Zt$ line bundle over a point), so that $x=1-N$.  Let
$L$ be the M\"obius line bundle on $\cs$, equipped with the $\Zt$-action 
that restricts to the trivial action over $\R^-$ and the
nontrivial action over $\R^+$, and put $e=1-L\in\Kcs$. The
equivariant Stiefel-Whitney classes of $N$ and $L$ generate
$H^1_{\Zt}\!(X;\Zt)$ \cite{Pr}, hence $x$ and $e$ generate $\LA$.
The relations $N^2 = L^2 = 1$ translate into $x^2=2x$ and $e^2 =
2e$. To obtain another relation, consider a pair of complexes
$$0\to 1\overset{g}\longrightarrow L\to 0\hspace{1cm}\text{ and }
\hspace{1cm}0\to N\overset{g'}\longrightarrow L\to 0$$
representing $e$ and $x-e$, respectively.
The map $g$ is forced to be zero over $\R^+$, but we may choose
it to be injective elsewhere.  Similarly, we may choose
$g'$ to vanish only on $\R^-$.
Tensoring these two complexes together, we obtain an exact complex
representing $e(x-e)$, hence this class is trivial in $\Kcs$.
In Theorem \ref{pl} we will prove that these are all of the relations.
\end{example}

Let $\etai = \oi^*e\in \LA$.  Equivariant line bundles on $\MA$ are
classified by the group $H^1_{\Zt}\!\big(\MA;\Zt\big)$, which is
generated by the pullbacks of the equivariant Stiefel-Whitney
classes of $L$ and $N$ along the various maps $\oi$ \cite{Pr}.  Then
by naturality of the equivariant Stiefel-Whitney class, $\LA$ is
generated multiplicitively by $\eta_1,\ldots,\eta_n$ and $x$.

\begin{remark}\label{conj}
We may rephrase this observation by saying that the pullback
$\omega^*:\Line\big((\cs)^n\big)\to\Line\big(\MA\big)$ along the map
$\omega = (\omega_1,\ldots,\omega_n):\MA\to (\cs)^n$ is surjective.
Note that if $\operatorname{rk}\A = d$, then $\omega$ is an embedding,
and $\omega^*$ is simply the restriction map.  
It seems reasonable to conjecture that the pullback
$\omega^*:\kog\big((\cs)^n\big)\to\kog\big(\MA\big)$
is surjective as well.
\end{remark}

For any connected component $C\subs\CA$,
let $h_C:\LA\to\Line(C) = \Z[x]/x(2-x)$ be the map given by restriction to $C$.

\begin{lemma}\label{restrict}
For all $C$, $h_C$ takes $\etai$ to $x$ if $C\subs H_i^+$, and to $0$
if $C\subs H_i^-$.
\end{lemma}

\begin{proof}
Restricting to the real locus commutes with pulling back along $\oi$,
hence it is enough to see that $e|_{\R^+}=x$ and $e|_{\R^-}=0$.
This observation follows from the representation of $e$ as a complex
in Example \ref{cstar}.
\end{proof}

\vspace{-\baselineskip}
\begin{definition}\label{pres}
Let $\PA$ be the ring
$\Z[e_1,\ldots,e_n,x]/\IA$, where
$\IA$ is generated by the following five families of
relations.\,\footnote{Note that all of these relations are polynomial;
the $x^{-1}$ in the fifth family of relations cancels with a factor of $x$.}
\begin{eqnarray*}
&1)&\hs\hs x(2-x)\\&&\\
&2)&\hs\hs e_i(2-e_i)\hs\hs\text{ for }\hs\hs i\in\{1,\ldots,n\}\\&&\\
&3)&\hs\hs e_i(e_i-x)\hs\hs\text{ for }\hs\hs i\in\{1,\ldots,n\}\\&&\\
&4)&\hs\hs
\prod_{i\in\Sp}e_i\times\prod_{j\in\Sm}(e_j-x)
\hs\hs\text{ if }\hs\hs \bigcap_{i\in\Sp}\Hip\cap\bigcap_{j\in\Sm}H_j^- = \varnothing\\&&\\
&5)&\hs\hs
x^{-1}\left(\prod_{i\in\Sp}e_i\times\prod_{j\in\Sm}(e_j-x)
-\prod_{i\in\Sp}(e_i-x)\times\prod_{j\in\Sm}e_j\right)
\text{ if }\hs\hs \bigcap_{i\in\Sp}\Hip\cap\bigcap_{j\in\Sm}H_j^-
=\varnothing\\
&&\hspace{15pt} \text{ and }\bigcap_{i\in S}H_i \text{ is nonempty
with codimension less than $|S|$, where $S = \Sp\sqcup\Sm$.}
\end{eqnarray*}
\end{definition}

\begin{remark}
For the fourth and fifth families of generators of $\IA$, it is sufficient
to consider only pairs of subsets $S^+,S^-\subs\on$ which are minimal
with respect to the given conditions; 
the other relations are generated by these.
\end{remark}

\begin{remark}\label{coor}
The notation that we have chosen is slightly abusive,
as the ideal $\IA$ and the ring $\PA$ depend not just on the
hyperplane arrangement, but also on the choice of linear
forms $\omega_i$ used to define the hyperplanes.  If $\omega_i$
is scaled by a positive real number, nothing changes, but if it is
scaled by a negative real number, the roles of $\Hip$
and $\Him$ are reversed.  Let $\A^i$ be the same arrangement
as $\A$ with the sign of $\omega_i$ reversed.  Then
$\PA$ is isomorphic to $P(\A^i)$ via the map $e_i\mapsto x-e_i$,
which justifies the abuse.  This point will be revisited in the beginning
of the proof of Theorem \ref{pl}.
\end{remark}

\begin{remark}
The generators of $\IA$ in Definition \ref{pres} are similar to the
relations given by Varchenko and Gelfand in their presentation
of the ring $\vgz$ \cite{VG}.  In fact, the ring $\PA/\!\<x-2\>$ is
isomorphic to the subring of $\vgz$ generated by two times the Heaviside
functions.
We give an abstract characterization of this subring for simple arrangements
in Section \ref{vg}.
\end{remark}

\begin{definition}
A {\em circuit} is a minimal set $S$ such that $\cap_{i\in S}H_i$
is nonempty with codimension less than $|S|$.
All circuits admit a unique decomposition $S=\Sp\sqcup\Sm$ (up to permutation of
the two pieces)
such that $$\bigcap_{i\in\Sp}\Hip\cap\bigcap_{j\in\Sm}H_j^-=\varnothing.$$
A set $T$ is called a {\em broken circuit} if there exists $i$
with $i<j$ for all $j\in T$ such that $T\cup\{i\}$ is a circuit.
An {\em nbc-set} $A\subs\on$ is a set such that $\cap_{i\in A}H_i$ is nonempty
and $A$ does not contain a broken circuit.
\end{definition}

For any subset $A\subs\on$, let $e_A = \prod_{i\in A}e_i$.

\begin{lemma}\label{free}
The ring $\PA$ is additively a free abelian group of rank $R+1$, where $R$
is the number of connected components of $\CA$.
\end{lemma}

\begin{proof}
The set $\{x\}\cup\{e_A\mid A \text{ an nbc-set}\}$ is an additive basis
for $\PA$.  The monomials indexed by nbc-sets also form a basis for the
Orlik-Solomon algebra $A(\A;\Z)$, which is free-abelian of rank $R$
(see for example \cite[\S 2]{Yu}).
Hence $\PA$ is free abelian of rank $R+1$.
\end{proof}

In the following theorem, we show that the relations between the
K-theory classes $\eta_1,\ldots,\eta_n,x\in\LA$ are exactly given
by the ideal $\IA$. 


\begin{theorem}\label{pl}
The homomorphism $\phi:\Z[e_1,\ldots,e_n,x]\to\LA$ given by $\phi(e_i)=\etai$
and $\phi(x)=x$ is surjective with kernel $\IA$, hence $\LA$ is
isomorphic to $\PA$.
\end{theorem}

\begin{proof}
We begin by reducing to the case in which the 
polyhedron  $\D=\cap_{i=1}^n \Him$ is nonempty.
To do this, let $\A^i$ be as in Remark \ref{coor}.  We then
obtain a diagram as follows,
\begin{figure}[h]
\centerline{
\xymatrix{ 
\PA\ar[rd]^{\Phi}
\ar[rr]
&& P(\A^i)\ar[dl]_{\Phi^i} \\ 
& \LA
}}
\end{figure}

\noindent
which commutes because $(-\omega_i)^*e = x-\eta_i$.
This tells us that $\Phi$ is an isomorphism if and only if $\Phi^i$
is an isomorphism.  By changing the signs of enough
of the linear forms, we may achieve the condition that $\D$
is nonempty.

To see that $\IA$ is contained in the kernel of $\phi$, we must show
that each of the families of generators maps to zero.
The images under $\phi$
of the first three families are all pullbacks of relations in $\Kcs$,
and are therefore zero in $\LA$.

Let $$Y_i^+ = \oi^{-1}(\R^+)\text{    and    }Y_i^- = \oi^{-1}(\R^-)\subs\MA.$$
We have already observed that $e\in\Kcs$ may be represented by a complex which
is exact away from $\R^+$, therefore $\etai = \oi^*(e)$ may be represented
by a complex which is exact away from $Y_i^+$.
Similarly, $\etai-x = \oi^*(e-x)$
may be represented by a complex which is exact away from $Y_i^-$.
Suppose that $$p\in\bigcap_{i\in\Sp}Y_i^+\cap\bigcap_{j\in\Sm}Y_j^-.$$
Then the real part $$\Re(p)\in\bigcap_{i\in\Sp}\Hip\cap\bigcap_{j\in\Sm}H_j^-,$$
hence
$$\bigcap_{i\in\Sp}\Hip\cap\bigcap_{j\in\Sm}H_j^-=\varnothing
\hspace{5pt}\impl\hspace{5pt}
\bigcap_{i\in\Sp}Y_i^+\cap\bigcap_{j\in\Sm}Y_j^-=\varnothing.$$
In this case $\prod_{i\in\Sp}\etai\times\prod_{j\in\Sm}(\etaj-x)$ is represented
by an exact complex, and is therefore equal to zero.
This accounts for the fourth family of generators of $\IA$.

Now suppose given a circuit $S = S^+\sqcup S^-\subs\on$ with
$\Big(\bigcap_{i\in\Sp}\Hip\Big)\cap\Big(\bigcap_{j\in\Sm}H_j^-\Big)=\varnothing$,
and consider the arrangement $\A_S = \{H_i^{\C}\mid i\in S\}$.
The space $$\MAS = \Cd\smallsetminus\bigcup_{i\in S}H_i$$ contains the space $\MA$,
and we have a commutative diagram
$$\begin{CD}
\Z[e_i,x]_{i\in S} @>\phi_S >> \LAS\\
@VVV @ VVV\\
\Z[e_1,\ldots,e_n,x] @>\phi >> \LA\\
\end{CD}$$
where the map from $\LAS$ to $\LA$ is given by restriction.  Hence to show that
the class
$$x^{-1}\left(\prod_{i\in\Sp}e_i\times\prod_{j\in\Sm}(e_j-x)
-\prod_{i\in\Sp}(e_i-x)\times\prod_{j\in\Sm}e_j\right)$$ is in the
kernel of $\phi$, it will suffice to show that it is in the kernel
of $\phi_S$. Dividing by the vector space $\bigcap_{i\in S}H_i$,
which is a factor of $\MAS$, we obtain a homotopy equivalent space
$\mathcal{M}(\hat\A_S)$, where $\hat\A_S$ is a central, essential
arrangement of $|S|$ hyperplanes in a vector space of dimension
$|S|-1$. Thus we have reduced to the special case where $\A$ is a
central arrangement of $d+1$ generic hyperplanes in $\R^d$.

Let us number our hyperplanes $H_0,\ldots,H_d$.
We have a splitting $\{0,\ldots,d\} = S^+\sqcup S^-$
such that $$\bigcap_{i\in\Sp}\Hip\cap\bigcap_{j\in\Sm}H_j^-=\varnothing,$$
and the fact that $\D = \bigcap_{i=0}^d H_i^-$ is nonempty implies that
either $S^+$ or $S^-$ is a singleton consisting of the unique hyperplane
that is not a facet of $\D$.  Without loss of generality, let us assume
that $S^-=\{0\}$.  Let $\dA$ be the decone of $\A$ with respect to $H_0$.
More explicitly, $\dA$ is the cooriented arrangement of $d$ affine hyperplanes in the
affine space $V = \{p\in\Cd\mid\omega_0(p)=-1\}$ whose hyperplanes are cut out
by the restrictions of $\omega_1,\ldots,\omega_d$ to $V$.
(We ask that $\omega_0(p)=-1$ rather than $1$ so that $\D\cap\Mda$ will
be nonempty.)

The ring $\Lda$ is generated
by $\nu_1,\ldots,\nu_d$ and $x$, where the generator $\nu_i$
corresponding to the hyperplane $H_i\cap V$ is equal to the restriction of $\etai$
to $\Mda\subs\MA$.  We have $$\bigcap_{i=1}^d (H_i^+\cap V) =
(-\D)\cap V = \varnothing,$$
hence we have the relation $\prod_{i=1}^d\nu_i = 0 \in\Lda$.

Consider the map $f:\MA\to\Mda$ given by the formula $f(p) = -p/\omega_0(p)$.

\begin{lemma}\label{fstar}
For all $i\in\{1,\ldots,d\}$, $f^*(\nu_i) = \eta_0 + \eta_i - \eta_0\eta_i$.
\end{lemma}

\begin{proof}
We will prove the equivalent statement that
$1-f^*(\nu_i)=(1-\etaz)(1-\etai)$. From the definitions of $x,
\nu_i$, and $\etai$, we see that both sides of this equation can be
represented by honest equivariant line bundles (rather than virtual
bundles), hence we may interpret the statement as an equation in the
Picard group of $\Zt$-equivariant line bundles on $\MA$.  This group
is isomorphic to the cohomology group $H^1_{\Zt}(\MA;\Zt)$, which
injects into $H^1_{\Zt}(\CA;\Zt)$ by the restriction map 
\cite[2.4 $\&$ 2.5]{Pr}.
Since the isomorphism between the Picard group and the first
equivariant cohomology commutes with restriction, it is enough to
prove that $$h_C\(1-f^*(\nu_i)\) = h_C\Big((1-\etaz)(1-\etai)\Big)$$
for all components $C\subs\CA$. By Lemma \ref{restrict}, and the
observation that $(1-x)^2 = 1$, we have
$$h_C\Big((1-\etaz)(1-\etai)\Big)=
\begin{cases}
1-x & \text{if $\omega_0$ and $\oi$ take values of opposite sign on $C$}\\
1 & \text{otherwise.}
\end{cases}$$
On the other hand, $f^*(\nu_i) = \(-\frac{\oi}{\omega_0}\)^*(e)$.
Using the fact that restriction to the real locus commutes with pulling back,
and the fact that $e|_{\R^+}=x$ and $e|_{\R^-}=0$, we obtain the desired equality.
\end{proof}

By Lemma \ref{fstar}, we have
\begin{equation}\label{first}
0 = \prod_{i=1}^d f^*(\nu_i) = \prod_{i=1}^d (\etaz + \etai -
\etaz\etai) = \prod_{i=1}^d \Big(\etai(1 - \etaz)+\etaz\Big)
=\sum_{A\subs\{1,\ldots,d\}}\!\! (1 -
\etaz)^{|A|}\cdot\etaz^{|A^c|}\cdot\eta_A,
\end{equation}
where $\eta_A = \prod_{i\in A}\etai$. Since $(1-\etaz)^2 = 1$ and
$(1-\etaz)\cdot \etaz = - \etaz$, we also have
\begin{equation}\label{second}
(1-\etaz)^{|A|}\cdot\etaz^{|A^c|} =
\begin{cases}
1-\etaz & \text{if $A=\{1,\ldots,d\}$ and $d$ is odd}\\
(-1)^{|A|}\cdot\etaz^{|A^c|} &\text{otherwise.}
\end{cases}
\end{equation}
On the other hand, consider the expression
$$x^{-1}\(\etaz\prod_{i=1}^d(\etai-x)-(\etaz-x)\prod_{i=1}^d\etai\),$$
which may be rewritten as
$$\sum_{A\subs\{1,\ldots,d\}}(-\etaz)^{|A^c|}\cdot\eta_A
\,=\, (-1)^d\cdot\!\!\!\sum_{A\subs\{1,\ldots,d\}}(-1)^{|A|}
\cdot\etaz^{|A^c|}\cdot\eta_A.$$ By Equations \eqref{first} and
\eqref{second}, this expression is equal to zero if $d$ is even, and
otherwise equal to
$$\prod_{i=1}^d\etai+(1-\etaz)\prod_{i=1}^d\etai
\,=\,(2-\etaz)\prod_{i=1}^d\etai \,=\,
(x-\etaz)\prod_{i=1}^d\etai.$$ But
$\displaystyle{H_0^-\cap\bigcap_{i=1}^d\Hip = \varnothing},$ therefore
$\displaystyle{(x-\etaz)\prod_{i=1}^d\etai=0}$ from the fourth
family of relations.  Hence we have shown that
$$x^{-1}\(\etaz\prod_{i=1}^d(\etai-x)-(\etaz-x)\prod_{i=1}^d\etai\)=0,$$
and therefore that all of the generators of $\IA$ are contained
in the kernel of $\phi$.

Our work up to this point implies that $\phi$ descends to a surjection
$\hat\phi:\PA\to\LA$;
it remains to show that $\hat\phi$ is injective.
We prove instead the following stronger statement.
Let $h:\LA\to\Lca$ be the restriction to the fixed point set.
(The ring $\Lca$ is a direct sum one copy of $\Z[x]/x(2-x)$ for each
component $C\subs\CA$, and $h$ is the direct sum of the maps $h_C$.)

\begin{lemma}\label{inject}
The composition $h\circ\hat\phi:\PA\to\Lca$ is injective.
\end{lemma}

\begin{proof}
By Lemma \ref{free}, it is enough to prove injectivity after tensoring
with the rational numbers $\Q$.
Given any component $C\subs\CA$, choose a pair of subsets $S^+,S^-\subs\on$
such that $\Big(\bigcap_{i\in\Sp}\Hip\Big)\cap\Big(\bigcap_{j\in\Sm}H_j^-\Big) = C$.
Then for any other component $D\subs\CA$, Lemma \ref{restrict} tells us that
$$h_D\(\prod_{i\in\Sp}\etai\cdot\prod_{j\in\Sm}(x-\etaj)\)
=\delta_{CD}\cdot x^{|S^+\cup S^-|},$$
hence $$h\(\prod_{i\in\Sp}\etai\cdot\prod_{j\in\Sm}(x-\etaj)\)$$
is supported on a single component of $\CA$.
The $R$ elements obtained this way, along with the trivial
vector bundle $1$, generate an $(R+1)$-dimensional subspace of
$\Line\big(\CA\big)\otimes\Q$.
Since $\dim\PA\otimes\Q = R+1$, $h\circ\hat\phi$ must be injective.
\end{proof}

Injectivity of $h\circ\hat\phi$ implies injectivity of $\hat\phi$,
therefore $\hat\phi:\PA\to\LA$ is an isomorphism.
This completes the proof of Theorem \ref{pl}.
\end{proof}
\end{section}

\begin{section}{The subring description of \boldmath$\LA$}\label{vg}

An arrangement $\A$ is called {\em simple}
if $\codim\cap_{i\in S}H_i = |S|$
for any $S$ such that $\cap_{i\in S}H_i$ is nonempty.
By Theorem \ref{pl} and Lemma \ref{inject}, we know
that $\PA\cong\LA$ is isomorphic to a subring of $\Line\big(\CA\big)$.
In this section we give a combinatorial interpretation of that subring
in the special case where the arrangement $\A$ is simple, which we will assume
for the rest of the section.

The arrangement $\A$ divides $\R^d$ into a polytopal complex $|\A|$
whose maximal faces are the connected components $C\subs\CA$,
and whose smaller faces are the open faces of the polytopes $\bar{C}$.
Given any face $F\in|\A|$, we let $\mathfrak{C}_F$ denote
the set of maximal faces $C$ containing $F$ in their closure,
and we choose a sign function $\eps_F:\mathfrak{C}_F\to\{\pm 1\}$
such that any two maximal faces separated by a single hyperplane
receive a different sign.

\begin{definition}
Let $\BA$ be the subgroup of $\Lca\cong\displaystyle{\bigoplus_{C\subs\CA}}\Z[x]/x(2-x)$
defined by the following condition.
$$\mu\in\BA\text{ if and only if for all faces }F\in|\A|,\hs
\sum_{C\in\mathfrak{C}_F}\eps_F(C)\mu_C\,\in\,\<x^{\codim F}\>\subs\Z[x]/x(2-x).$$
\end{definition}

It is clear that $\BA$ is a subgroup.  It is not so obvious that it
is also a subring; this fact will follow from Theorem \ref{pb}.

\begin{proposition}\label{lands}
The image of the restriction map $h:\LA\to\Lca$ is contained in $\BA$.
\end{proposition}

\begin{proof}
We need only check that $h(\eta_A)\in\BA$ for all subsets
$A\subs\on$.  Choose a face $F\in|A|$, and a component $C\in\mathfrak{C}_F$.
Lemma \ref{restrict} tells us that
$$h_C(\eta_A)=
\begin{cases}
x^{|A|} & \text{ if $C\subs\bigcap_{i\in A}$}\Hip\\
0 & \text{ otherwise.}
\end{cases}$$
If $|A|\geq\codim F$, we are done.  If not, then by simplicity,
there must be an index $j$ such that $F\subs H_j$ but $j\notin A$.
In this case, $h_C(\eta_A) = h_D(\eta_A)$, where
$D\in\mathfrak{C}_F$ is the component separated from $C$ by $H_j$.
Since $\eps_F(C) = -\eps_F(D)$, these two terms of
$\displaystyle{\sum_{C\in\mathfrak{C}_F}}\eps_F(C)h_C(\eta_A)$
will cancel with each other.  Thus every term will cancel the contribution
of another term, and the total sum will be zero.
\end{proof}

Given a cooriented arrangement $\A=\{H_1,\ldots,H_n\}$
in $\R^d$, let $$\A' = \{H_1,\ldots,H_{n-1}\}$$ denote the arrangement
obtained by deleting $H_n$, and let
$$\A'' = \{H_i\cap H_n\mid i<n \text{ and }H_i\cap H_n\neq\varnothing\}$$
denote the arrangement of hyperplanes in $H_n$ given by restriction.
If $H_n\cap\D\neq \varnothing$, then $\A'$ and $\A''$ remain in the class of arrangements
with $\D$ nonempty.

\begin{proposition}\label{exact}
We have an exact sequence of groups
$$0\to P(\A')\overset{\a}{\longrightarrow} \PA\overset{\b}{\longrightarrow}
P(\A'')\overset{\g}{\longrightarrow}\Z\to 0.$$
The map $\a$ is the ring homomorphism taking $e_i$ to $e_i$ and $x$ to $x$.
We define $\b$ on the additive basis
$\{x\}\cup\{e_A\mid A \text{ an nbc-set}\}$
by $\b(x)=0$ and
$$\b(e_A) =
\begin{cases}
e_{A\smallsetminus\{n\}} & \text{if $n\in A$}\\
0 & \text{otherwise.}
\end{cases}$$
The third map $\g$ is defined by extracting the coefficient of $x$
in the corresponding basis for $P(\A'')$.
\end{proposition}

\begin{proof}
Injectivity of $\a$ is a consequence of the fact that every nbc-set for $\A'$
is also an nbc-set for $\A$.
Similarly, exactness at $\PA$ and $P(\A'')$ follow from the fact that if $n\in A$,
then  $A$ is an nbc-set for $\A$ if and only if $A\smallsetminus\{n\}$ is an nbc-set
for $\A''$.  Surjectivity of $\g$ is trivial.
\end{proof}

\vspace{-\baselineskip}
\begin{remark}
Our proof of Proposition \ref{exact} holds for arbitrary arrangements,
not just simple ones.
\end{remark}

Suppose that $H_n\cap\D$ is nonempty, and consider the sequence
$$0\to B(\A')\overset{a}{\longrightarrow} \BA\overset{b}{\longrightarrow}
B(\A'')\overset{c}{\longrightarrow}\Z\to 0$$
defined as follows.
The map $a$ is given by restriction, and $c$ is given
by taking the coefficient of $x$ corresponding to the component
$\D\cap H^n$ of $\Capr$.
Given an element $\mu\in\BA$ and a component $C''$ of $\Capr$,
we put $b(\mu)|_{C''} = (\mu_C-\mu_D)/x$, where $C\subs H_n^+$ and
$D\subs H_n^-$ are the two components of $\CA$ neighboring $C''$.
The fact that $\mu_C-\mu_D$ is a multiple of $x$ follows from the fact
that $\mu\in\BA$.  There is an inherent ambiguity in dividing by $x$,
owing to the fact that $x$ is annihilated by $2-x$.  We resolve this ambiguity
by requiring that $b(\mu)|_{\D\cap H^n}\in\Z$, and $b(\mu)|_{C''}$
is congruent to $b(\mu)|_{\D\cap H^n}$ modulo $x$ for all components $C''$.
This sequence is evidently a complex, and it is easy to check exactness
at $B(\A')$, $\BA$, and $\Z$.  Exactness at $B(\A'')$ will fall out of the process
of proving the following theorem.

\begin{proposition}\label{pb}
The restriction map $h:\LA\to\BA$ is an isomorphism.
\end{proposition}

\begin{proof}
By Theorem \ref{pl}, it is sufficient to prove that the composition
$h\circ\hat\phi:\PA\to\BA$ is an isomorphism.
We proceed by induction on the number of hyperplanes.  The base case
$n=0$ is trivial.
For the inductive step,
%
consider the commutative diagram
$$\xymatrix{
0 \ar[r] & P(\A) \ar[r]^{\a}\ar[d] & \PA \ar[r]^{\b} \ar[d]
& P(\A'') \ar[r]^{\g}\ar[d] & \Z \ar[r]\ar[d] & 0\\
0 \ar[r] & B(\A') \ar[r]^a & \BA \ar[r]^b & B(\A'')\ar[r]^c & \Z \ar[r] &0,}
$$
where the first three downward arrows are given by the composition
$h\circ\hat\phi$, and the last is the identity map.
Our inductive hypothesis tells us that the maps from $P(\A')$ to $B(\A')$
and $P(\A'')$ to $B(\A'')$ are isomorphisms.  This, along with
exactness of the top row at $P(\A'')$, implies the exactness
of the bottom row at $B(\A'')$.  Our Theorem then follows from the Five Lemma.
\end{proof}
\end{section}

\end{spacing}
\end{document}